\documentclass{article}

\usepackage[english]{babel}    
\usepackage{amsmath}           
\usepackage{amssymb}           
\usepackage[utf8]{inputenc}    
\usepackage[T1]{fontenc}       
\usepackage{eucal}             
\usepackage{bbm}               
\usepackage{longtable}         
\usepackage{exscale}           
\usepackage{mathtools}         
\usepackage[sort]{cite}        
\usepackage{mathrsfs}          
\usepackage{array}             
\usepackage{stmaryrd}          
\usepackage{paralist}          
\usepackage{wasysym}           
\usepackage{fancyhdr}          
\usepackage{xspace}            

 \newtheorem{thm}{Theorem}[section]

 \newtheorem{defn}[thm]{Definition}

 \newtheorem{ex}{Example}
 \numberwithin{equation}{section}

\newcommand{\I}          {\mathrm{i}}
\newcommand{\E}          {\mathrm{e}}
\newcommand{\D}          {\operatorname{\mathrm{d}}}
\newcommand{\cc}[1]      {\overline{{#1}}}
\newcommand{\id}         {\operatorname{\mathsf{id}}}

\newcommand{\Lie}        {\operatorname{\mathscr{L}\!}}

\newcommand{\Pol}        {\operatorname{\mathrm{Pol}}}

\newcommand{\End}        {\operatorname{\mathsf{End}}}
\newcommand{\SP}[1]      {\left\langle{#1}\right\rangle}
\newcommand{\ring}[1]    {\mathsf{#1}}

\newcommand{\lie}[1]     {\mathfrak{#1}}

\newcommand{\at}[1]      {\big|_{#1}}
\newcommand{\At}[1]      {\Big|_{#1}}

\newcommand{\Anti}       {\Lambda}
\newcommand{\Sym}        {\mathrm{S}}

\newcommand{\Schouten}[1]{\left\llbracket{#1}\right\rrbracket}

\newcommand{\starweyl}  {\mathbin{\star_{\scriptscriptstyle\mathrm{Weyl}}}}
\newcommand{\starstd}   {\mathbin{\star_{\scriptscriptstyle\mathrm{Std}}}}

\newcommand{\weylrep}  {\varrho_{\scriptscriptstyle\mathrm{Weyl}}}
\newcommand{\stdrep}   {\varrho_{\scriptscriptstyle\mathrm{Std}}}

\newcommand{\Sec}[1][k]      {\Gamma^{#1}}
\newcommand{\Secinfty}       {\Sec[\infty]}
\newcommand{\Fun}[1][k]      {\mathscr{C}^{#1}}
\newcommand{\Cinfty}         {\Fun[\infty]}

\newcommand{\argument}       {\,\cdot\,}

\newcommand{\tensor}[1][{}]           {\mathbin{\otimes_{\scriptscriptstyle{#1}}}}
\newcommand{\algebra}[1]      {\mathcal{#1}}
\newcommand{\module}[1]       {\mathcal{#1}}
\newcommand{\bimodule}[5]    {\sideset{^{\scriptscriptstyle{#1}}_{\scriptscriptstyle{#2}}}{^{\scriptscriptstyle{#4}}_{\scriptscriptstyle{#5}}}{\operatorname{#3}}}
\newcommand{\BEA}            {\bimodule{}{\algebra{B}}{\module{E}}{}{\algebra{A}}}
\newcommand{\AEpB}           {\bimodule{}{\algebra{A}}{\module{E}}{\prime}{\algebra{B}}}
\newcommand{\EA}             {\bimodule{}{}{\module{E}}{}{\algebra{A}}}

\begin{document}

%
%
%
%
%
%
%
%
%

\title{Recent Developments in Deformation Quantization}

\author{Stefan Waldmann\thanks{stefan.waldmann@mathematik.uni-wuerzburg.de}\\[0.5cm]
  Institut für Mathematik \\
  Lehrstuhl für Mathematik X \\
  Universität Würzburg \\
  Campus Hubland Nord \\
  Emil-Fischer-Straße 31 \\
  97074 Würzburg \\
  Germany}



\date{January 2015}
\maketitle
\begin{abstract}
    In this review an overview on some recent developments in
    deformation quantization is given.  After a general historical
    overview we motivate the basic definitions of star products and
    their equivalences both from a mathematical and a physical point
    of view. Then we focus on two topics: the Morita classification of
    star product algebras and convergence issues which lead to the
    nuclear Weyl algebra.
\end{abstract}

\section{Introduction: a Historical Tour D'Horizon}

In the last decades, deformation quantization evolved into a widely
accepted quantization scheme which, on one hand, provides deep
conceptual insights into the question of quantization and, on the
other hand, proved to be a reliably technique leading to explicit
understanding of many examples. It will be the aim of this review to
give some overview on the developments of deformation quantization
starting from the beginnings but also including some more recent
ideas.

The original formulations of deformation quantization by Bayen
et. al. aimed mainly at finite-dimensional classical mechanical
systems described by symplectic or Poisson manifolds
\cite{bayen.et.al:1978a} and axiomatize the heuristic quantization
formulas found earlier by Weyl, Groenewold, and Moyal \cite{weyl:1931a,
  groenewold:1946a, moyal:1949a}. Berezin considered the more
particular case of bounded domains and Kähler manifolds
\cite{berezin:1975c, berezin:1975b, berezin:1975a}. Shortly after it
proved to be a valuable tool to approach also problems in quantum
field theories, see e.g. the early works of Dito \cite{dito:1990a,
  dito:1992a, dito:1993a}.

Meanwhile, the question of existence and classification of deformation
quantizations, i.e. of star products, on symplectic manifolds was
settled: first DeWilde and Lecomte showed the existence of star
products on symplectic manifolds \cite{dewilde.lecomte:1983b} in 1983
after more particular classes \cite{dewilde.lecomte:1983a,
  dewilde.lecomte:1983c} were considered. Remarkably, also in 1983 the
first genuine class of Poisson structures was shown to admit star
products, the linear Poisson structures on the dual of a Lie algebra,
by Gutt \cite{gutt:1983a} and Drinfel'd \cite{drinfeld:1983a}. In 1986
Fedosov gave a very explicit and constructive way to obtain star
products on a symplectic manifold by means of a symplectic connection
\cite{fedosov:1986a}, see also \cite{fedosov:1996a, fedosov:1994a} for
a more detailed version. His construction is still one of the
cornerstones in deformation quantization as it provides not only a
particularly nice construction allowing to adjust many special
features of star products depending on the underlying manifold like
e.g. separation of variables (Wick type) on Kähler manifolds
\cite{karabegov:2000a, karabegov:1996a, bordemann.waldmann:1997a,
  neumaier:2003a} or star products on cotangent bundles
\cite{bordemann.neumaier.waldmann:1998a,
  bordemann.neumaier.waldmann:1999a,
  bordemann.neumaier.pflaum.waldmann:2003a}. Even beyond the
symplectic world, Fedosov's construction was used to globalize the
existence proofs of star products on Poisson manifolds
\cite{cattaneo.felder.tomassini:2002b, dolgushev:2005a}.  Yet another
proof for the symplectic case was given
in\cite{omori.maeda.yoshioka:1991a}.

Even though the symplectic case was understood well, the question of
existence on Poisson manifolds kept its secrets till the advent of
Kontsevich's formality theorem, solving his formality conjecture
\cite{kontsevich:1997a, kontsevich:1997:pre, kontsevich:2003a}. To
give a adequate overview on Kontsevich's formality theorem would
clearly go beyond the scope of this short review. Here one can rely on
various other publications like e.g. \cite{cattaneo:2005a,
  esposito:2015a}. In a nutshell, the formality theorem proves a very
general fact about smooth functions on a manifold from which it
follows that every (formal series of) Poisson structures can be
quantized into a star product, including a classification of star
products. Parallel to Kontsevich's groundbreaking result, the
classification of star products on symplectic manifolds was achieved
and compared by several groups \cite{bertelson.cahen.gutt:1997a,
  nest.tsygan:1995a, deligne:1995a, gutt.rawnsley:1999a,
  neumaier:2002a}. Shortly after Kontsevich, Tamarkin gave yet another
approach to the quantization problem on Poisson manifolds
\cite{tamarkin:1998b}, see also \cite{kontsevich:1999a,
  kontsevich.soibelman:2000a}, based on the language of operads and
the usage of Drinfel'd associators.  Starting with these formulations,
formality theory has evolved and entered large areas of contemporary
mathematics, see e.g. \cite{alekseev.meinrenken:2006a,
  alekseev.torossain:2010a, alekseev.torossian:2012a,
  dolgushev.tamarkin.tsygan:2007a,dolgushev.rubtsov:2009a,
  kontsevich:2001a, kontsevich:1999a} to name just a few.

While deformation quantization undoubtedly gave many important
contribution to pure mathematics over the last decades, it is now
increasingly used in contemporary quantum physics as well: perhaps
starting with the works of Dütsch and Fredenhagen on the perturbative
formulations of algebraic quantum field theory
\cite{duetsch.fredenhagen:2001b, duetsch.fredenhagen:1999a,
  duetsch.fredenhagen:2001a} it became clear that star products
provide the right tool to formulate quantum field theories in a
semiclassical way, i.e. as formal power series in $\hbar$. Now this
has been done in increasing generalities for various scenarios
including field theories on general globally hyperbolic spacetimes,
see e.g. \cite{brunetti.fredenhagen.verch:2003a,
  brunetti.fredenhagen.ribeiro:2012a:pre, baer.ginoux.pfaeffle:2007a,
  hollands.wald:2010a}.

Of course, from a physical point of view, deformation quantization can
not yet be the final answer as one always deals with formal power
series in the deformation parameter $\hbar$. A physically reasonable
quantum theory, however, requires of course convergence. Again, in the
very early works \cite{bayen.et.al:1978a} some special cases were
treated, namely the Weyl-Moyal product for which an integral formula
exists which allows for a reasonable analysis based on the Schwartz
space. The aims here are at least two-fold. On one hand one wants to
establish a reasonable spectral calculus for particular elements in
the star product algebra which allows to compute spectra in a
physically sensitive way. This can be done with the star exponential
formalism, which works in particular examples but lacks a general
framework. On the other hand, one can try to establish from the formal
star product a convergent version such that in the end one obtains a
$C^*$-algebra of quantum observables being a deformation, now in a
continuous way, of the classical functions on the phase space. This is
the point of view taken by strict deformation quantization, most
notably advocated by Rieffel \cite{rieffel:1989a, rieffel:1993a} and
Landsman \cite{landsman:1998a}, see also
\cite{bordemann.meinrenken.schlichenmaier:1991a,
  cahen.gutt.rawnsley:1995a, cahen.gutt.rawnsley:1994a,
  cahen.gutt.rawnsley:1993a, cahen.gutt.rawnsley:1990a} for the
particular case of quantizable Kähler manifolds and
\cite{natsume.nest.peter:2003a, natsume.nest:1999a, natsume:2000a} for
more general symplectic manifolds. Bieliavsky and coworkers found a
generalization of Rieffel's approach by passing from actions of the
abelian group $\mathbb{R}^d$ to more general Lie group actions
\cite{bieliavsky.detournay.spindel:2009a, bieliavsky:2002a,
  bieliavsky.massar:2001b}.  Having a $C^*$-algebra one has then the
full power of $C^*$-algebra techniques at hands which easily allows to
get a reasonable spectral calculus. However, constructing
$C^*$-algebraic quantizations is still very much in development: here
one has not yet a clear picture on the existence and classification of
the quantizations. In fact, one even has several competing definitions
of what one is looking for. It is one of the ongoing research projects
by several groups to understand the transition between formal and
strict quantizations in more detail.

Needless to say, that in the above historical survey we can barely
scratch on the surface of this vast topic: many aspects have not been
mentioned like the role played by symmetries and reduction, the
applications to concrete physical systems, various generalizations of
deformation quantization to other geometric brackets, alternative
approaches to various index theorems, relations to noncommutative
geometry, and many more.  In the remaining part of this review we will
focus on two aspects of the theory: first, we discuss the role of
classification results beyond the notion of equivalence,
i.e. isomorphism. Here we are particularly interested in the
classification of star products up to Morita equivalence. Second, we
give a short outlook on star products in infinite dimensions and
problems arising there by investigating one particular example: the
Weyl algebra of a vector space with a (quite arbitrary) bilinear
form. Beside the purely algebraic construction we obtain a locally
convex algebraic deformation once we start in this category.

%
%

\section{From Poisson Manifolds to Star Products}

In this section we give a more detailed but still non-technical
motivation of the definition of star products and list some first
examples.

The set-up will be a finite-dimensional phase space which we model by
a symplectic or, more generally, a Poisson manifold $(M, \pi)$ where
$\pi \in \Secinfty(\Anti^2 TM)$ is a bivector field satisfying
\begin{equation}
    \label{eq:Jacobi}
    \Schouten{\pi, \pi} = 0.
\end{equation}
Here $\Schouten{\argument, \argument}$ is the Schouten bracket and the
condition is equivalent to the Jacobi identity for the Poisson bracket
\begin{equation}
    \label{eq:PoissonBracket}
    \{f, g\} = - \Schouten{\Schouten{f, \pi}, g} = \pi(\D f, \D g)
\end{equation}
determined by $\pi$ for functions $f, g \in \Cinfty(M)$. One can then
formulate classical Hamiltonian mechanics using $\pi$ and
$\{\argument, \argument\}$. For a gentle introduction to Poisson
geometry see \cite{waldmann:2007a} as well as \cite{vaisman:1994a,
  dufour.zung:2005a, laurent-gengoux.pichereau.vanhaecke:2013a,
  cannasdasilva.weinstein:1999a}. One has several important examples
of Poisson manifolds:
\begin{itemize}
\item Every symplectic manifold $(M, \omega)$ where $\omega \in
    \Secinfty(\Anti^2T^*M)$ is a closed non-degenerate two-form, is a
    Poisson manifold with $\pi = \omega^{-1}$. The Jacobi identity
    \eqref{eq:Jacobi} corresponds then directly to $\D\omega = 0$.
\item Every cotangent bundle $T^*Q$ is a symplectic manifold in a
    canonical way with an exact symplectic form $\omega = \D \theta$
    where $\theta \in \Secinfty(T^*(T^*Q))$ is the canonical (or
    tautological) one-form on $T^*Q$. This is the arena of geometric
    mechanics.
\item Kähler manifolds are particularly nice examples of symplectic
    manifolds as they possess a compatible Riemannian metric and a
    compatible complex structure.
\item The dual $\lie{g}^*$ of a Lie algebra $\lie{g}$ is always a
    Poisson manifold with a linear Poisson structure: the coefficient
    functions of the tensor field $\pi$ are linear functions on
    $\lie{g}^*$, explicitly given by
    \begin{equation}
        \label{eq:KKS}
        \{f, g\}(x)
        =
        x_i c^i_{k\ell}
        \frac{\partial f}{\partial x_k}
        \frac{\partial g}{\partial x_\ell},
    \end{equation}
    where $x_1, \ldots, x_n$ are the linear coordinates on $\lie{g}^*$
    and $c^i_{k\ell}$ are the corresponding structure constants of
    $\lie{g}$. Since it vanishes at the origin, this is never
    symplectic.
\item Remarkably and slightly less trivial is the observation that on
    every manifold $M$ and every $p \in M$ there is a Poisson
    structure $\pi$ with compact support but $\pi\at{p}$ has maximal
    rank.
\end{itemize}

To motivate the definition of a star product we consider the most easy
example of the classical phase space $\mathbb{R}^2$ with canonical
coordinates $(q, p)$. Canonical quantization says that we have to map
the spacial coordinate $q$ to the position operator $Q$ acting on a
suitable domain in $L^2(\mathbb{R}, \D x)$ as multiplication
operator. Moreover, we have to assign the momentum coordinate $p$ to
the momentum operator $P = - \I \hbar \frac{\partial}{\partial q}$,
again defined on a suitable domain. Since we want to ignore
functional-analytic questions at the moment, we simply chose
$\Cinfty_0(\mathbb{R})$ as common domain for both operators. In a next
step we want to quantize polynomials in $q$ and $p$ as well. Here we
face the ordering problem as $pq = qp$ but $PQ \ne QP$. One simple
choice is the \emph{standard ordering}
\begin{equation}
    \label{eq:StandardOrd}
    q^np^m \; \mapsto \;
    \stdrep(q^np^m) = Q^nP^m
    =
    (-\I\hbar)^m q^n\frac{\partial^m}{\partial q^m}
\end{equation}
for monomials and its linear extension to all polynomials. More
explicitly, this gives
\begin{equation}
    \label{eq:StdRep}
    \stdrep(f) = \sum_{r=0}^\infty \frac{1}{r!}
    \left(\frac{\hbar}{\I}\right)^r
    \frac{\partial^r f}{\partial p^r}\At{p=0}
    \frac{\partial^r}{\partial q^r}.
\end{equation}
Now this formula still makes sense for smooth functions $f$ which are
polynomial in $p$, i.e. for $f \in \Cinfty(\mathbb{R})[p]$.  The main
idea of deformation quantization is now to pull-back the operator
product: this is possible since the image of $\stdrep$ is the space of
all differential operators with smooth coefficients which therefore is
a (noncommutative) algebra. We define the \emph{standard-ordered star
  product} by
\begin{equation}
    \label{eq:starstd}
    f \starstd g
    = \stdrep^{-1}(\stdrep(f)\stdrep(g))
    = \sum_{r=0}^\infty \frac{1}{r!}
    \left(\frac{\hbar}{\I}\right)^r
    \frac{\partial^r f}{\partial p^r}
    \frac{\partial^r g}{\partial q^r}
\end{equation}
for $f, g \in \Cinfty(\mathbb{R})[p]$. While it is clear that
$\starstd$ is an associative product the behaviour with respect to the
complex conjugation is bad: we do not get a $^*$-involution,
$\cc{f \starstd g} \ne \cc{g} \starstd \cc{f}$, since
\begin{equation}
    \label{eq:NeumaierOp}
    \stdrep(f)^\dag = \stdrep(N^2 f)
    \quad
    \textrm{with}
    \quad
    N
    =
    \exp\left(
        \frac{\hbar}{2\I} \frac{\partial^2}{\partial q \partial p}
    \right),
\end{equation}
as a simple integration by parts shows. We can repair this unpleasant
feature by defining the Weyl ordering and the Weyl product by
\begin{equation}
    \label{eq:WeylStuff}
    \weylrep(f) = \stdrep(Nf)
    \quad
    \textrm{and}
    \quad
    f \starweyl g = N^{-1}(Nf \starstd Ng).
\end{equation}
Note that $N$ is indeed an invertible operator on
$\Cinfty(\mathbb{R})[p]$. Again, $\starweyl$ is associative. Then we
get
\begin{equation}
    \label{eq:WeylFeatures}
    \cc{f \starweyl g} = \cc{g} \starweyl \cc{f}
    \quad
    \textrm{and}
    \quad
    \weylrep(f \starweyl g) = \weylrep(f) \weylrep(g).
\end{equation}
For both products we can collect the terms of order $\hbar^r$ which
gives
\begin{equation}
    \label{eq:StarProd}
    f \star g = \sum_{r=0}^\infty \hbar^r C_r(f, g)
\end{equation}
with bidifferential operators $C_r$ of order $r$ in each argument. The
explicit formula for $\starweyl$ is slightly more complicated than the
one for $\starstd$ in \eqref{eq:starstd} but still easy to compute. We
have
\begin{equation}
    \label{eq:StarFirstOrders}
    f \star g = fg + \cdots
    \quad
    \textrm{and}
    \quad
    f \star g - g \star f = \I\hbar \{f, g\} + \cdots,
\end{equation}
where $+ \cdots$ means higher orders in $\hbar$. Also $f \star 1 = f =
1 \star f$. Note that the seemingly infinite series in
\eqref{eq:StarProd} is always finite as long as we take functions in
$\Cinfty(\mathbb{R})[p]$.

The idea is now to axiomatize these features for $\star$ in such a way
that it makes sense to speak of a star product on a general Poisson
manifold. The first obstacle is that on a generic manifold $M$ there
is nothing like functions which are polynomial in certain
coordinates. This is a chart-dependent characterization which one does
not want to use. But then already for $\starweyl$ and $\starstd$ one
encounters the problem that for general $f, g \in
\Cinfty(\mathbb{R}^2)$ the formulas \eqref{eq:starstd} and
\eqref{eq:WeylFeatures} will not make any sense: the series are indeed
infinite and since we can adjust the Taylor coefficients of a smooth
function in a rather nasty way, there is no hope for convergence.
The way out is to consider \emph{formal} star product in a first step,
i.e. formal power series in $\hbar$. This yields the definition of
star products \cite{bayen.et.al:1978a}:
\begin{defn}
    A formal star product $\star$ on a Poisson manifold $(M, \pi)$ is
    an associative $\mathbb{C}[[\hbar]]$-bilinear associative product
    for  $\Cinfty(M)[[\hbar]]$ such that
    \begin{equation}
        \label{eq:FormalStarProd}
        f \star g = \sum_{r=0}^\infty \hbar^r C_r(f, g)
    \end{equation}
    with
    \begin{enumerate}
    \item $C_0(f, g) = fg$,
    \item $C_1(f, g) - C_1(g, f) = \I\{f, g\}$,
    \item $C_r(1, f) = 0 = C_r(f, 1)$ for $r \ge 1$,
    \item $C_r$ is a bidifferential operator.
    \end{enumerate}
\end{defn}
Already in the trivial example above we have see that there might be
more than one star product. The operator $N$ interpolates between them
and is invisible in classical physics: for $\hbar = 0$ the operator
$N$ becomes the identity. As a formal series of differential operator
starting with the identity it is invertible and implements an algebra
isomorphism. This is now taken as definition for equivalence of star
products: given two star products $\star$ and $\star'$ on a manifold,
a formal power series $T = \id + \sum_{r=1}^\infty \hbar^r T_r$ of
differential operators $T_r$ with $T1 = 1$ is called an equivalence
between $\star$ and $\star'$ if one has
\begin{equation}
    \label{eq:Equivalence}
    f \star' g = T^{-1}(Tf \star Tg).
\end{equation}
Note that $T$ is indeed invertible as a formal power series. Hence
this is an equivalence relation. Conversely, given such a $T$ and
$\star$ we get a new star product $\star'$ by \eqref{eq:Equivalence}.

After the general set-up we are now in the position to list some basic
examples of star products:
\begin{itemize}
\item The explicit formulas for $\starstd$ and $\starweyl$ immediately
    generalize to higher dimensions yielding equivalent star products
    on $\mathbb{R}^{2n}$ and hence also on every open subset of
    $\mathbb{R}^{2n}$. Since by the Darboux Theorem every symplectic
    manifold looks like an open subset of $\mathbb{R}^{2n}$
    \emph{locally}, the question of existence of star products on
    symplectic manifolds is a global problem.
\item For the linear Poisson structure \eqref{eq:KKS} on the dual
    $\lie{g}^*$ of a Lie algebra $\lie{g}$ one gets a star product as
    follows \cite{gutt:1983a}: First, we note that
    $\Sym^\bullet(\lie{g}) = \Pol^\bullet(\lie{g}^*)$. Then the
    PBW isomorphism
    \begin{equation}
        \label{eq:Symmetrizer}
        \Sym^\bullet(\lie{g}) \ni \xi_1 \vee \cdots \vee \xi_k
        \; \mapsto \;
        \frac{(\I\hbar)^k}{k!} \sum_{\sigma \in S_k}
        \xi_{\sigma(1)} \cdots \xi_{\sigma(k)}
        \in \mathcal{U}(\lie{g})
    \end{equation}
    from the symmetric algebra over $\lie{g}$ into the universal
    enveloping algebra allows to pull the product of
    $\mathcal{U}(\lie{g})$ back to $\Sym^\bullet(\lie{g})$ and hence
    to polynomials on $\lie{g}^*$. One can now show that after
    interpreting $\hbar$ as a formal parameter one obtains indeed a
    star product quantizing the linear Poisson bracket. This star
    product is completely characterized by the feature that
    \begin{equation}
        \label{eq:BCH}
        \exp(\hbar\xi) \star \exp(\hbar\eta)
        =
        \exp(\mathrm{BCH}(\hbar\xi, \hbar\eta))
    \end{equation}
    for $\xi, \eta \in \lie{g}$ with the Baker-Campbell-Hausdorff
    series $\mathrm{BCH}$, see \cite{gutt:1983a,
      bordemann.neumaier.waldmann:1998a}.
\item The next interesting example is perhaps the complex projective
    space $\mathbb{CP}^n$ and its non-compact dual, the Poincaré disc
    $\mathbb{D}_n$ with their canonical Kähler structures of constant
    holomorphic sectional curvature. For these, star products were
    considered by Moreno and Ortega-Navarro
    \cite{moreno.ortega-navarro:1983b} who gave recursive formulas
    using local coordinates. Cahen, Gutt, and Rawnsley
    \cite{cahen.gutt.rawnsley:1995a, cahen.gutt.rawnsley:1994a,
      cahen.gutt.rawnsley:1993a, cahen.gutt.rawnsley:1990a} discussed
    this in their series of papers of quantization of Kähler manifolds
    as one of the examples. The first explicit (non-recursive) formula
    was found in \cite{bordemann.brischle.emmrich.waldmann:1996a,
      bordemann.brischle.emmrich.waldmann:1996b} by a quantization of
    phase space reduction and extended to complex Grassmannians in
    \cite{schirmer:1997a:pre}. Ever since these star products have
    been re-discovered by various authors.
\end{itemize}

We briefly comment on the general existence results: as already
mentioned, the symplectic case was settled in the early 80s. The
Poisson case follows from Kontsevich's formality theorem.
\begin{thm}[Kontsevich]
    On every Poisson manifold there exist star products.
\end{thm}
The classification is slightly more difficult to describe: we consider
\emph{formal Poisson structures}
\begin{equation}
    \label{eq:FormalPoisson}
    \pi = \hbar \pi_1 + \hbar^2 \pi_2 + \cdots
    \in \hbar \Secinfty(\Anti^2 TM)[[\hbar]]
    \quad
    \textrm{with}
    \quad
    \Schouten{\pi, \pi} = 0.
\end{equation}
Moreover, let $X = \hbar X_1 + \hbar^2 X_2 + \cdots \in
\hbar\Secinfty(TM)[[\hbar]]$ be a formal vector field, starting in
first order of $\hbar$. Then one calls $\exp(\Lie_X)$ a formal
diffeomorphism which defines an action
\begin{equation}
    \label{eq:expLX}
    \exp(\Lie_X)\colon
    \Secinfty(\Anti^2 TM)[[\hbar]] \ni \nu
    \; \mapsto \;
    \nu + \Lie_X\nu + \frac{1}{2} \Lie_X^2 \nu + \cdots
    \in
    \Secinfty(\Anti^2 TM)[[\hbar]].
\end{equation}
Via the Baker-Campbell-Hausdorff series, the set of formal
diffeomorphisms becomes a group and \eqref{eq:expLX} is a group
action. Since $\Lie_X$ is a derivation of the Schouten bracket, it
follows that the action of $\exp(\Lie_X)$ preserves formal Poisson
structures. The space of orbits of formal Poisson structures modulo
this group action gives now the classification:
\begin{thm}[Kontsevich]
    The set of equivalence classes of formal star products is in
    bijection to the set of equivalence classes of formal Poisson
    structures modulo formal diffeomorphisms.
\end{thm}
In general, both moduli spaces are extremely difficult to describe.
However, if the first order term $\pi_1$ in $\pi$ is symplectic, then
we have a much easier description which is in fact entirely
topological:
\begin{thm}[Bertelson, Cahen, Gutt, Nest, Tsygan, Deligne, \ldots]
    On a symplectic manifold $(M, \omega)$ the equivalence classes of
    star products are in bijection to the formal series in the second
    deRham cohomology. In fact, one has a canonical surjective map
    \begin{equation}
        \label{eq:CharClass}
        c\colon
        \star \; \mapsto \;
        c(\star) \in
        \frac{[\omega]}{\I\hbar} +
        \mathrm{H}^2_{\mathrm{dR}}(M, \mathbb{C})[[\hbar]]
    \end{equation}
    such that $\star$ and $\star'$ are equivalent iff $c(\star) =
    c(\star')$.
\end{thm}
This map is now called the characteristic class of the symplectic star
product. In a sense which can be made very precise
\cite{bursztyn.dolgushev.waldmann:2012a}, the inverse of $c(\star)$
corresponds to Kontsevich's classification by formal Poisson tensors.

%
%

\section{Morita Classification}
\label{sec:Morita}

We come now to some more particular topics in deformation
quantization. In this section we discuss a coarser classification
result than the above classification up to equivalence.

The physical motivation to look for Morita theory is rather simple and
obvious: in quantum theory we can not solely rely on the observable
algebra as the only object of interest. Instead, we also need to have
a reasonable notion of states. While for $C^*$-algebras there is a
simple definition of a state as a normalized positive functional, in
deformation quantization we do not have $C^*$-algebras in a first
step. Surprisingly, the notion of positive functionals still makes
sense if interpreted in the sense of the ring-ordering of
$\mathbb{R}[[\hbar]]$ and it produces a physically reasonable
definition, see \cite{bordemann.waldmann:1998a}. However, the
requirements from quantum theory do not stop here: we also need a
super-position principle for states. Since positive functionals can
only be added convexly, we need to realize the positive functionals as
expectation value functionals for a $^*$-representation of the
observable algebra on some (pre) Hilbert space. Then we can take
complex linear combination of the corresponding vectors to implement
the super-position principle. This leads to the need to understand the
representation theory of the star product algebras, a program which
was investigated in great detail \cite{bursztyn.waldmann:2001b,
  bursztyn.waldmann:2001a, bursztyn.waldmann:2004a,
  bursztyn.waldmann:2005b, bursztyn.dolgushev.waldmann:2012a,
  jansen.waldmann:2006a}, see also \cite{waldmann:2005b} for a review.
The main point is that replacing the ring of scalars from $\mathbb{R}$
to $\mathbb{R}[[\hbar]]$ and thus from $\mathbb{C}$ to
$\mathbb{C}[[\hbar]]$ works surprisingly well as long as we do not try
to implement analytic concepts: the non-archimedean order of
$\mathbb{R}[[\hbar]]$ forbids a reasonable analysis. However, the
concept of positivity is entirely algebraic and hence can be used and
employed in this framework as well.

In fact, one needs not to stop here: \emph{any} ordered ring
$\ring{R}$ instead of $\mathbb{R}$ will do the job and one can study
$^*$-algebras over $\ring{C} = \ring{R}(\I)$ and their
$^*$-representation theory on pre Hilbert modules over $\ring{C}$.
For many reasons it will also be advantageous to consider
representation spaces where the inner product is not taking values in
the scalars but in some \emph{auxiliary} $^*$-algebra $\mathcal{D}$.
\begin{ex}
    Let $E \longrightarrow M$ be a complex vector bundle over a smooth
    manifold $M$. Then $\Secinfty(E)$ is a $\Cinfty(M)$-module in the
    usual way. A Hermitian fiber metric $h$ give now a sesquilinear
    map
    \begin{equation}
        \label{eq:FiberMetric}
        \SP{\argument, \argument}\colon
        \Secinfty(E) \times \Secinfty(E) \longrightarrow \Cinfty(M)
    \end{equation}
    which is also $\Cinfty(M)$-linear in the second argument, i.e. we
    have $\SP{s, tf} = \SP{s, t}f$ for all $s, t \in \Secinfty(E)$ and
    $f \in \Cinfty(M)$. Moreover, the pointwise positivity of $h_p$ on
    $E_p$ implies that the map
    \begin{equation}
        \label{eq:CompletePositive}
        \SP{\argument, \argument}^{(n)}\colon
        \Secinfty(E)^n \times \Secinfty(E)^n
        \longrightarrow M_n(\Cinfty(M)) = \Cinfty(M, M_n(\mathbb{C}))
    \end{equation}
    is positive for all $n$ in the sense that the matrix-valued
    function $\SP{S, S}^{(n)} \in \Cinfty(M, M_n(\mathbb{C}))$ yields
    a positive matrix at all points of $M$ for all $S = (s_1, \ldots,
    s_n) \in \Secinfty(E)^n$.
\end{ex}

Using this kind of \emph{complete positivity} for an inner product
yields the definition of a pre Hilbert right module over a
$^*$-algebra $\mathcal{D}$, where the inner product takes values in
$\mathcal{D}$. Then again, we can formulate what are
$^*$-representations of a $^*$-algebra $\mathcal{A}$ on such a pre
Hilbert right module over $\mathcal{D}$. Without further difficulties
this gives various categories of $^*$-representations of $^*$-algebras
on inner product modules or pre Hilbert modules over auxiliary
$^*$-algebras.

Having now a good notion of $^*$-representations of $^*$-algebras it
is a major task to understand the resulting categories for those
$^*$-algebras occurring in deformation quantization. Now from
$C^*$-algebra theory we anticipate that already with the full power of
functional-analytic techniques it will in general be impossible to
``understand'' the category of $^*$-representations completely, beside
rather trivial examples. The reason is that there will simply be too
many inequivalent such $^*$-representations and a decomposition theory
into irreducible ones is typically an extremely hard problem. In a
purely algebraic situation like for formal star product algebras,
things are even worse: here we expect even more inequivalent ones
which are just artifacts of the algebraic formulation. There are many
examples of inequivalent $^*$-representations which, after one
implements mild notions of convergence and hence of analytic aspects,
become equivalent. From a physical point of view such inequivalences
would then be negligible. However, it seems to be quite difficult to
decide this \emph{before} convergence is implemented, i.e. on the
algebraic side.

Is the whole program now useless, hopeless? The surprising news is
that one can indeed say something non-trivial about the
$^*$-representation theories of the star product algebras from
deformation quantization, and for $^*$-algebras in general. The idea
is that even if the $^*$-representation theory of a given $^*$-algebra
is horribly complicated and contains maybe unwanted
$^*$-representation, we can still \emph{compare} the whole
$^*$-representation theory of one $^*$-algebra to another $^*$-algebra
and ask whether they are equivalent as categories.

This is now the basic task of Morita theory. To get a first impression
we neglect the additional structure of ordered rings,
$^*$-involutions, and positivity and consider just associative
algebras over a common ring of scalars. For two such algebras
$\algebra{A}$ and $\algebra{B}$ we want to know whether their
categories of left modules are equivalent categories. Now there might
be many very strange functors implementing an equivalence and hence
one requires them to be compatible with direct sums of modules, which
is clearly a reasonable assumption. The prototype of such a functor is
then given by the tensor product with a $(\algebra{B},
\algebra{A})$-bimodule. Since the tensor product with $\algebra{A}$
itself is (for unital algebras) naturally isomorphic to the identity
functor and since the tensor product of bimodules is associative up to
a natural isomorphism, the question of equivalence of categories via
such tensor product functors becomes equivalent to the question of
\emph{invertible bimodules}: Here a $(\algebra{B},
\algebra{A})$-bimodule $\BEA$ is called invertible if there is a
$(\algebra{A}, \algebra{B})$-bimodule $\AEpB$ such that the tensor
product $\BEA \tensor[\algebra{A}] \AEpB$ is isomorphic to
$\algebra{B}$ and $\AEpB \tensor[\algebra{B}] \BEA$ is isomorphic to
$\algebra{B}$, always as bimodules.

The classical theorem of Morita now gives a complete and fairly easy
description of the possible bimodules with this property: $\BEA$ has
to be a finitely generated projective and full right
$\algebra{A}$-module and $\algebra{B}$ is isomorphic to
$\End_{\algebra{A}}(\EA)$ via the left module structure, see
e.g. \cite{lam:1999a}.

Now the question is how such bimodules look like for star product
algebras. Classically, the finitely generated projective modules over
$\Cinfty(M)$ are, up to isomorphism, just sections $\Secinfty(E)$ of a
vector bundle $E \longrightarrow M$, this is the famous Serre-Swan
theorem in its incarnation for differential geometry. As soon as the
fiber dimension is non-zero, the fullness condition is trivially
satisfied. Hence the only Morita equivalent algebras to $\Cinfty(M)$
are, again up to isomorphism, the sections $\Secinfty(\End(E))$ of
endomorphism bundles. The corresponding bimodule is then
$\Secinfty(E)$ on which both algebras act in the usual way. It
requires now a little argument to see that for star products, an
equivalence bimodule gives an equivalence bimodule in the classical
limit $\hbar = 0$, i.e. a vector bundle. Conversely, the sections of
every vector bundle can be deformed into a right module over the star
product algebra in a unique way up to isomorphism. Thus for star
products, we have to look for the corresponding module endomorphisms
of such deformed sections of vector bundles. Finally, in order to get
again a star product algebra, the endomorphisms of the deformed
sections have to be, in the classical limit, isomorphic to the
functions on a manifold again. This can only happen if the vector
bundle was actually a line bundle over the same manifold. Hence the
remaining task is to actually compute the star product of the algebra
acting from the left side when the star product for the algebra on the
right side is known. Here one has the following results:
\begin{thm}[Bursztyn, W. \cite{bursztyn.waldmann:2002a}]
    Let $(M, \omega)$ and $(M', \omega')$ be a symplectic manifolds
    and let $\star$, $\star'$ be two star products on $M$ and $M'$,
    respectively. Then $\star$ and $\star'$ are Morita equivalent iff
    there exists a symplectomorphism $\psi\colon M \longrightarrow M'$
    such that
    \begin{equation}
        \label{eq:ClassesME}
        \psi^*c(\star') - c(\star) \in 2\pi\I
        \mathrm{H}^2_{\mathrm{dR}}(M, \mathbb{Z}).
    \end{equation}
    The difference of the above classes defines then a line bundle
    which implements the Morita equivalence bimodule by deforming its
    sections.
\end{thm}

This theorem has already an important physical interpretation: for
cotangent bundles $T^*Q$ the characteristic classes $c(\star)$ can be
interpreted as the classes of magnetic fields $B$ on the configuration
space $Q$. Then a quantization of a charged particle in the background
field of such a $B$ requires a star product with characteristic class
$c(\star)$. Compared to the trivial characteristic class, $c(\star) =
0$, the above theorem then tells that quantization with magnetic field
has the same representation theory iff the magnetic field satisfies
the integrality condition for a Dirac monopole. Thus we get a Morita
theoretic interpretation of the charge quantization for magnetic
monopoles which is now extremely robust against details of the
quantization procedure: the statement holds for all cotangent bundles
and for all equivalent star products with the given characteristic
class.

Also in the more general Poisson case the full classification is
known. Here the actual statement is slightly more technical as it
requires the Kontsevich class of the star products and a canonically
given action of the deRham cohomology on equivalence classes of formal
Poisson structures by gauge transformations. Then one obtains the
following statement, see also \cite{jurco.schupp.wess:2002a} for an
earlier heuristic argument based on noncommutative field theories:
\begin{thm}[Bursztyn, Dolgushev, W. \cite{bursztyn.dolgushev.waldmann:2012a}]
    Star products on Poisson manifolds are Morita equivalent iff their
    Kontsevich classes of formal Poisson tensors are gauge equivalent
    by a $2\pi\I$-integral deRham class.
\end{thm}

%
%

\section{Beyond Formal Star Products}
\label{sec:Beyond}

Since formal star products are clearly not sufficient for physical
purposes, one has to go beyond formal power series. Here several
options are available: on one hand one can replace the formal series
in the star products by integral formulas. The formal series can then
be seen as the asymptotic expansions of the integral formulas in the
sense of Taylor series of smooth functions of $\hbar$, which are
typically not analytic: hence we can not expect convergence.
Nevertheless, the integral formulas allow for a good analytic
framework.

However, if one moves to field theories and hence to
infinite-dimensional systems, quantization becomes much more
complicated. Surprisingly, series formulas for star products can still
make sense in certain examples, quite unlike the integral formulas:
such integrals would consist of integrations over a
infinite-dimensional phase space. Hence we know that such things can
hardly exist in a mathematically sound way.

This motivates the second alternative, namely to investigate the
formal series in the star products directly without integral formulas
in the back. This might also be possible in infinite dimensions and
yields reasonable quantizations there. While this is a program far
from being understood, we can present here now one class of examples
with a particular physical relevance: the Weyl algebra.

Here we consider a real vector space $V$ with a bilinear map
$\Lambda\colon V \times V \longrightarrow \mathbb{C}$. Then we
consider the complexified symmetric algebra
$\Sym_{\mathbb{C}}^\bullet(V)$ of $V$ and interpret this as the
polynomials on the dual $V^*$. In finite dimensions this is correct,
in infinite dimensions the symmetric algebra is better to be
interpreted as the polynomials on the (not necessarily existing)
pre-dual. On $V^*$, there are simply much more polynomials than the
ones arising from $\Sym^\bullet_{\mathbb{C}}(V)$. Now we can extend
$\Lambda$ to a biderivation
\begin{equation}
    \label{eq:PLambda}
    P_\Lambda\colon
    \Sym_{\mathbb{C}}^\bullet(V)
    \tensor
    \Sym_{\mathbb{C}}^\bullet(V)
    \longrightarrow
    \Sym_{\mathbb{C}}^\bullet(V)
    \tensor
    \Sym_{\mathbb{C}}^\bullet(V)
\end{equation}
in a unique way by enforcing the Leibniz rule in both tensor
factors. If we denote by $\mu\colon \Sym_{\mathbb{C}}^\bullet(V)
\tensor \Sym_{\mathbb{C}}^\bullet(V) \longrightarrow
\Sym_{\mathbb{C}}^\bullet(V)$ the symmetric tensor product, then
\begin{equation}
    \label{eq:PoissonLambda}
    \{a, b\}_\Lambda
    =
    \mu \circ (P_\Lambda (a \tensor b) - P_\Lambda(b \tensor a))
\end{equation}
is a Poisson bracket. In fact, this is the unique constant Poisson
bracket with the property that for linear elements $v, w \in V$ we
have $\{v, w\} = \Lambda(v, w) - \Lambda(w, v)$. Hence the
antisymmetric part of $\Lambda$ determines the bracket. However, we
will use the symmetric part for defining the star product. This will
allow to include also standard-orderings or other orderings like Wick
ordering from the beginning.

A star product quantizing this constant Poisson structure can then
be found easily. We set
\begin{equation}
    \label{eq:Star}
    a \star b = \mu \circ \exp(z P_\Lambda) (a \tensor b)
\end{equation}
where $z \in \mathbb{C}$ is the deformation parameter. For physical
applications we will have to set $z = \I\hbar$ later on.  Note that
$\star$ is indeed well-defined since on elements in the symmetric
algebra, the operator $P_\Lambda$ lowers the degree by one in each
tensor factor.

In a next step we want to extend this product to more interesting
functions than the polynomial-like ones. The strategy is to look for a
topology which makes the product continuous and which allows for a
large completion of $\Sym^\bullet_{\mathbb{C}}(V)$. To start with, one
has to assume that $V$ is endowed with a topology itself. Hence let
$V$ be a locally convex Hausdorff space. In typical examples from
quantum mechanics, $V$ is the (dual of the) phase space and hence
finite dimensional, which makes the topology unique. In quantum field
theory, $V$ would be something like a test function space, i.e. either
the Schwartz space $\mathcal{S}(\mathbb{R}^d)$ or $\Cinfty_0(M)$ for a
manifold $M$, etc. In this case $V$ would be a Fréchet or LF space.

We use now the continuous seminorms of $V$ to extend them to tensor
powers $V^{\tensor k}$ for all $k \in \mathbb{N}$ by taking their
tensor powers: we equip $V^{\tensor k}$ with the $\pi$-topology
inherited from $V$. This means that for a continuous seminorm $p$ on
$V$ we consider $p^{\tensor k}$ on $V^{\tensor k}$ and take all such
seminorms to define a locally convex topology on $V^{\tensor
  k}$. Viewing the symmetric tensor powers as a subspace, this induces
the $\pi$-topology also for $\Sym^\bullet_{\mathbb{C}}(V)$, simply by
restricting the seminorms $p^{\tensor k}$. For the whole symmetric
algebra we need to extend these seminorms we have on each symmetric
degree. This can be done in many inequivalent ways. Useful for our
purposes is the following construction. We fix a parameter $R \ge
\frac{1}{2}$ and define
\begin{equation}
    \label{eq:SeminormS}
    p_R(a) = \sum_{k=0}^\infty k!^R p^{\tensor k}(a_k)
\end{equation}
for every $a = \sum_{k=0}^\infty a_k$ with $a_k \in
\Sym^k_{\mathbb{C}}(V)$. Note that the sum is finite as long as we
take $a$ in the symmetric algebra. Now taking all those seminorms
$p_R$ for all continuous seminorms $p$ of $V$ induces a locally convex
topology on $V$. Clearly, this is again Hausdorff. Moreover, all
$\Sym^k_{\mathbb{C}}(V)$ are closed embedded subspaces in
$\Sym^\bullet_{\mathbb{C}}(V)$ with respect to this topology.

The remarkable property of this topology is now that a continuous
$\Lambda$ will induce a continuous star product \cite{waldmann:2014a}:
\begin{thm}
    Let $\Lambda\colon V \times V \longrightarrow \mathbb{C}$ be a
    continuous bilinear form on $V$. Then $\star$ is a continuous
    associative product on $\Sym^\bullet_{\mathbb{C}}(V)$ with respect
    to the locally convex topology induced by all the seminorms $p_R$
    with $p$ being a continuous seminorm on $V$, as long as $R \ge
    \frac{1}{2}$.
\end{thm}
The proof consists in an explicit estimate for $p_R(a \star b)$. Note
that the topology can \emph{not} be locally multiplicatively convex
since in the Weyl algebra we have elements satisfying canonical
commutation relations, thereby forbidding a submultiplicative
seminorm.
\begin{defn}[Locally convex Weyl algebra]
    Let $\Lambda\colon V \times V \longrightarrow \mathbb{C}$ be a
    continuous bilinear form on $V$. Then the completion of
    $\Sym^\bullet_{\mathbb{C}}(V)$ with respect to the above locally
    convex topology and with the canonical extension of $\star$ is
    called the locally convex Weyl algebra $\mathcal{W}_R(V, \star)$.
\end{defn}
Thus we have found a framework where the Weyl star product actually
converges. Without proofs we list a few properties of this Weyl
algebra:
\begin{itemize}
\item The locally convex Weyl algebra $\mathcal{W}_R(V, \star)$ is a
    locally convex unital associative algebra. The product $a \star b$
    can be written as the absolutely convergent series
    \begin{equation}
        \label{eq:SeriesWeyl}
        a \star b = \mu \circ \exp(z P_\Lambda) (a \tensor b).
    \end{equation}
\item The product $\star$ depends holomorphically on $z \in
    \mathbb{C}$.
\item For $\frac{1}{2} \le R < 1$ the locally convex Weyl algebra
    $\mathcal{W}_R(V, \star)$ contains the exponential functions
    $\E^{\alpha v}$ for all $v \in V$ and all $\alpha \in
    \mathbb{C}$. They satisfy the usual Weyl relations. Note that not
    only the unitary ones, i.e. for $\alpha$ imaginary, are contained
    in the Weyl algebra, but all exponentials.
\item The locally convex Weyl algebra is nuclear iff $V$ is
    nuclear. In all relevant examples in quantum theory this will be
    the case. In this case we refer to the \emph{nuclear Weyl
      algebra}.
\item If $V$ admits an absolute Schauder basis, then the symmetrized
    tensor products of the basis vectors constitute an absolute
    Schauder basis for the Weyl algebra, too.  Again, in many
    situations $V$ has such a basis.
\item The Weyl algebras for different $\Lambda$ on $V$ are isomorphic
    if the antisymmetric parts of the bilinear forms coincide.
\item Evaluations at points in the topological dual $V'$ are
    continuous linear functionals on $\mathcal{W}_R(V, \star)$. Hence
    we still can view the elements of the completion as particular
    functions on $V'$.
\item The translations by elements in $V'$ still act on
    $\mathcal{W}_R(V, \star)$ by continuous automorphisms. If $R < 1$
    these translations are inner automorphism as soon as the element
    $\varphi \in V'$ is in the image of the musical map induced by
    $\Lambda$.
\end{itemize}

We conclude this section now with a few comments on examples. First it
is clear that in finite dimensions we can take $V = \mathbb{R}^{2n}$
with the canonical Poisson bracket on the symmetric algebra. Then many
types of orderings can be incorporated in fixing the symmetric part of
$\Lambda$, while the antisymmetric part is given by the Poisson
bracket. Thus all the resulting star products allow for this analytic
framework. This includes examples known earlier in the literature, see
e.g. \cite{omori.maeda.miyazaki.yoshioka:2007a,
  beiser.roemer.waldmann:2007a}. In this case we get a nuclear Weyl
algebra with an absolute Schauder basis.

More interesting is of course the infinite dimensional case. Here we
have to specify the space $V$ and the bilinear form $\Lambda$ more
carefully. In fact, the \emph{continuity} of $\Lambda$ becomes now a
strong conditions since bilinear maps in locally convex analysis tend
to be only separately continuous without being continuous. However,
there are several situations where we can either conclude the
continuity of a bilinear separately continuous map by abstract
arguments, like for Fréchet spaces. Or one can show directly that the
particular bilinear form one is interested in is continuous. We give
one of the most relevant examples for (quantum) field theory:
\begin{ex}
    Let $M$ be a globally hyperbolic spacetime and let $D$ be a
    normally hyperbolic differential operator acting on a real vector
    bundle $E$ with fiber metric $h$. Moreover, we assume that $D$ is
    a connection Laplacian for a metric connection with respect to $h$
    plus some symmetric operator $B$ of order zero. In all relevant
    examples this is easy to obtain. Then one has advanced and
    retarded Green operators leading to the propagator $F_M$ acting on
    test sections $\Secinfty_0(E^*)$. We take $V = \Secinfty_0(E^*)$
    with its usual LF topology. Then
    \begin{equation}
        \label{eq:LambdaCov}
        \Lambda(\varphi, \psi)
        =
        \int_M h^{-1}(F_M(\varphi),  \psi) \mu_g
    \end{equation}
    is the bilinear form leading to the Peierls bracket on the
    symmetric algebra $\Sym^\bullet(V)$. Here $\mu_g$ is the metric
    density as usual. The kernel theorem then guarantees that
    $\Lambda$ is continuous as needed. Thus we obtain a locally convex
    and in fact nuclear Weyl algebra from this. Now $\Lambda$ is
    highly degenerated. It follows that in the Poisson algebra there
    are many Casimir elements. The kernel of $F_M$ generates a Poisson
    ideal and also an ideal in the Weyl algebra, which coincides with
    the vanishing ideal of the solution space. Hence dividing by this
    (Poisson) ideal gives a Poisson algebra or Weyl algebra which can
    be interpreted as the observables of the (quantum) field theory
    determined by the wave equation $Du = 0$. It can then be shown
    that for every Cauchy surface $\Sigma$ in $M$ there is a canonical
    algebra isomorphism to the Weyl algebra build from the symplectic
    Poisson algebra on the initial conditions on $\Sigma$. Details of
    this construction can be found in \cite{waldmann:2014a}, see also
    \cite{baer.ginoux.pfaeffle:2007a} for the background information
    on the wave equation.
\end{ex}


\subsection*{Acknowledgment}

It is a pleasure to thank the organizers of the Regensburg conference
for their kind invitation and the fantastic organization of this
stimulating conference. Moreover, I would like to that Chiara Esposito
for helpful remarks on the manuscript.


\end{document}